\documentclass[11pt,twoside]{amsart}

\usepackage{amssymb, palatino,graphicx, hyperref}
\usepackage[normalem]{ulem}
\usepackage[all]{xy}
\usepackage{xcolor}
\usepackage{ifthen}

\graphicspath{{./}{./Figures/}}

\def\F{\ensuremath{\mathcal{F}}}

\def\P{\ensuremath{\mathbb P}}

\def\Q{\ensuremath{\mathbb Q}}
\def\R{\ensuremath{\mathbb R}}
\def\Z{\ensuremath{\mathbb Z}}

\def\C{\ensuremath{\mathbb C}}

\def\F{\ensuremath{\mathbb F}}

\def\CP1{\ensuremath{\mathbb C \mathbb P^1}}

\def\Pn-1{\ensuremath{\P^{n-1}}}

\definecolor{red}{rgb}{.6,0,0}
\definecolor{green}{rgb}{0,.6,0}
\definecolor{darkgreen}{rgb}{0,0.3,0}
\definecolor{purple}{rgb}{0.5,0,0.5}
\definecolor{darkblue}{rgb}{0,0,0.7}
\definecolor{greenblue}{rgb}{0,0.4,0.5}
\definecolor{myblue}{HTML}{1685d6}
\definecolor{mypurple}{HTML}{b82e6c}

\newboolean{draft}

\newcommand{\cmt}[1]
{\ifthenelse {\boolean{draft}}
{{\sc \tiny \color{red} #1}}
{}}

\newcommand{\newbb}[1]
{\ifthenelse {\boolean{draft}}
{{\color{darkblue} #1}}
{#1}}

\newcommand{\newbbb}[1]
{\ifthenelse {\boolean{draft}}
{{\color{greenblue} #1}}
{#1}}

\newcommand{\nopost}[1]
{\ifthenelse {\boolean{draft}}
{{\color{cyan} #1}}
{}}

\newcommand{\maynopost}[1]
{\ifthenelse {\boolean{draft}}
{{\color{purple} #1}}
{}}

\newcommand{\margincmt}[1]
{\ifthenelse {\boolean{draft}}
{\marginpar{{\sc \tiny \color{red} #1}}}
{}}
\newcommand{\inred}[1]
{\ifthenelse{\boolean{draft}}{{\color{red} #1}}{#1}}

\newcommand{\new}[1]
{\ifthenelse {\boolean{draft}}
{{\color{green} #1}}
{#1}}

\newcommand{\neww}[1]
{\ifthenelse {\boolean{draft}}
{{\color{darkgreen} #1}}
{#1}}

\newcommand{\newb}[1]
{\ifthenelse {\boolean{draft}}
{{\color{blue} #1}}
{#1}}

\newcommand{\del}[1]
{\ifthenelse {\boolean{draft}}
{{\color{magenta} #1}}
{}}

\newboolean{details_on}

\newcommand{\details}[1]
{\ifthenelse {\boolean{details_on}}
{{\color{darkgreen} \tiny #1}}
{}}

\title{Hirzebruch surfaces in a one--parameter family}
\author{Fiammetta Battaglia, Elisa Prato, Dan Zaffran}
\thanks{}

\newcommand\vz{\ensuremath{\underline{z}}}

\newcommand{\bT}{\ensuremath{\mathcal{T} } } 

\newcommand{\calF}{\ensuremath{ \mathcal{F} } }

\setboolean{draft}{true}
\setboolean{details_on}{true}

\begin{document}
\maketitle
\begin{abstract} 
We introduce a family of spaces, parametrized by positive real numbers, that includes all of the Hirzebruch surfaces. 
Each space is viewed from two distinct perspectives. First, as a leaf space of a compact, complex, foliated manifold, following 
\cite{bz1}. Second, as a symplectic cut of the manifold $\C\times S^2$ in a possibly nonrational direction, following \cite{cut}. 
\end{abstract}
\section*{Introduction}
This article is dedicated to the memory of Paolo de Bartolomeis and presents a theme in which complex and symplectic geometry 
are closely intertwined. 

\begin{figure}[h]
\begin{center}
\includegraphics[width=7cm]{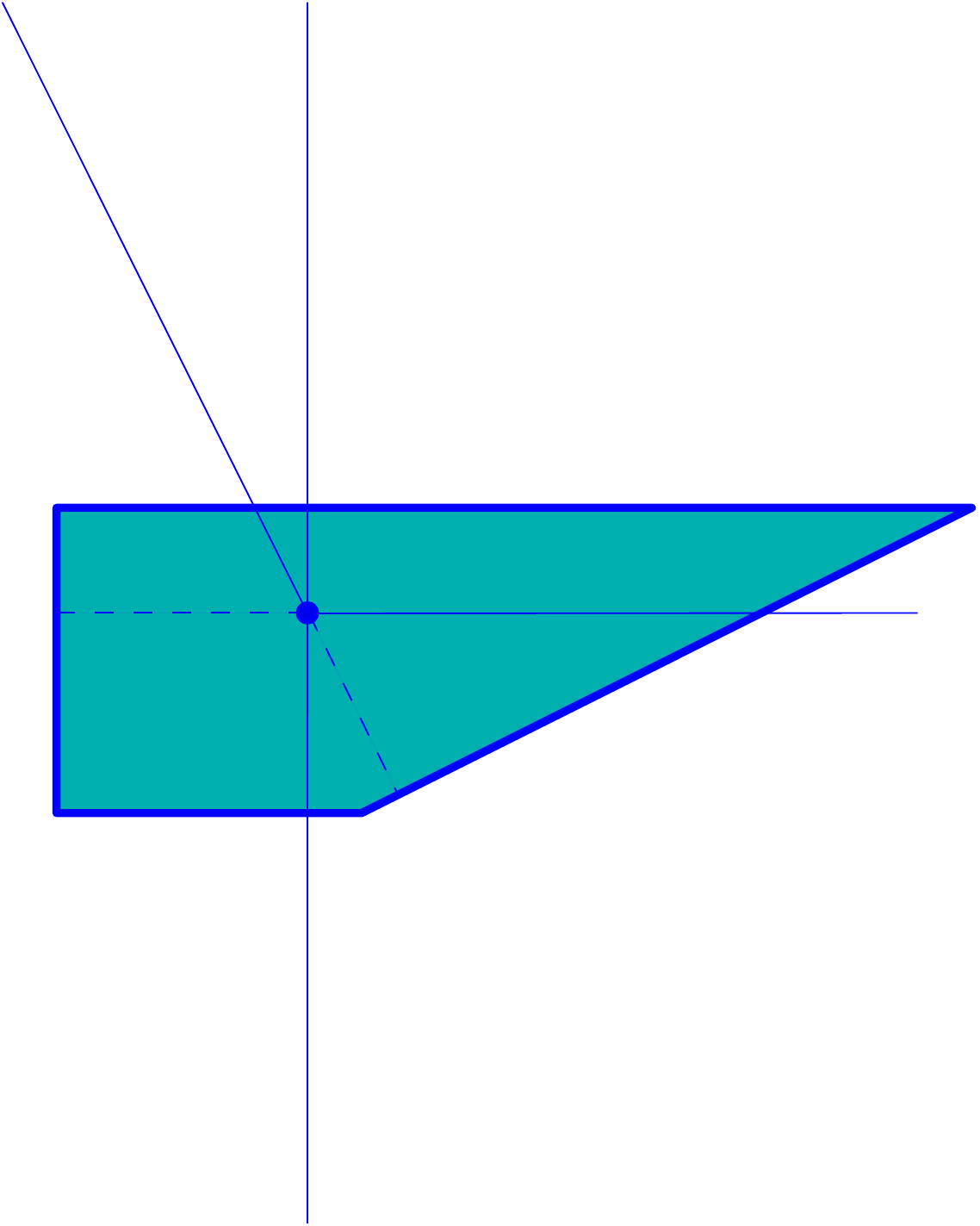}
\caption{The polytope $P_2$ with its normal fan.}
\label{smoothpolytope}
\end{center}
\end{figure}

\begin{figure}[h]
\begin{center}
\includegraphics[width=7cm]{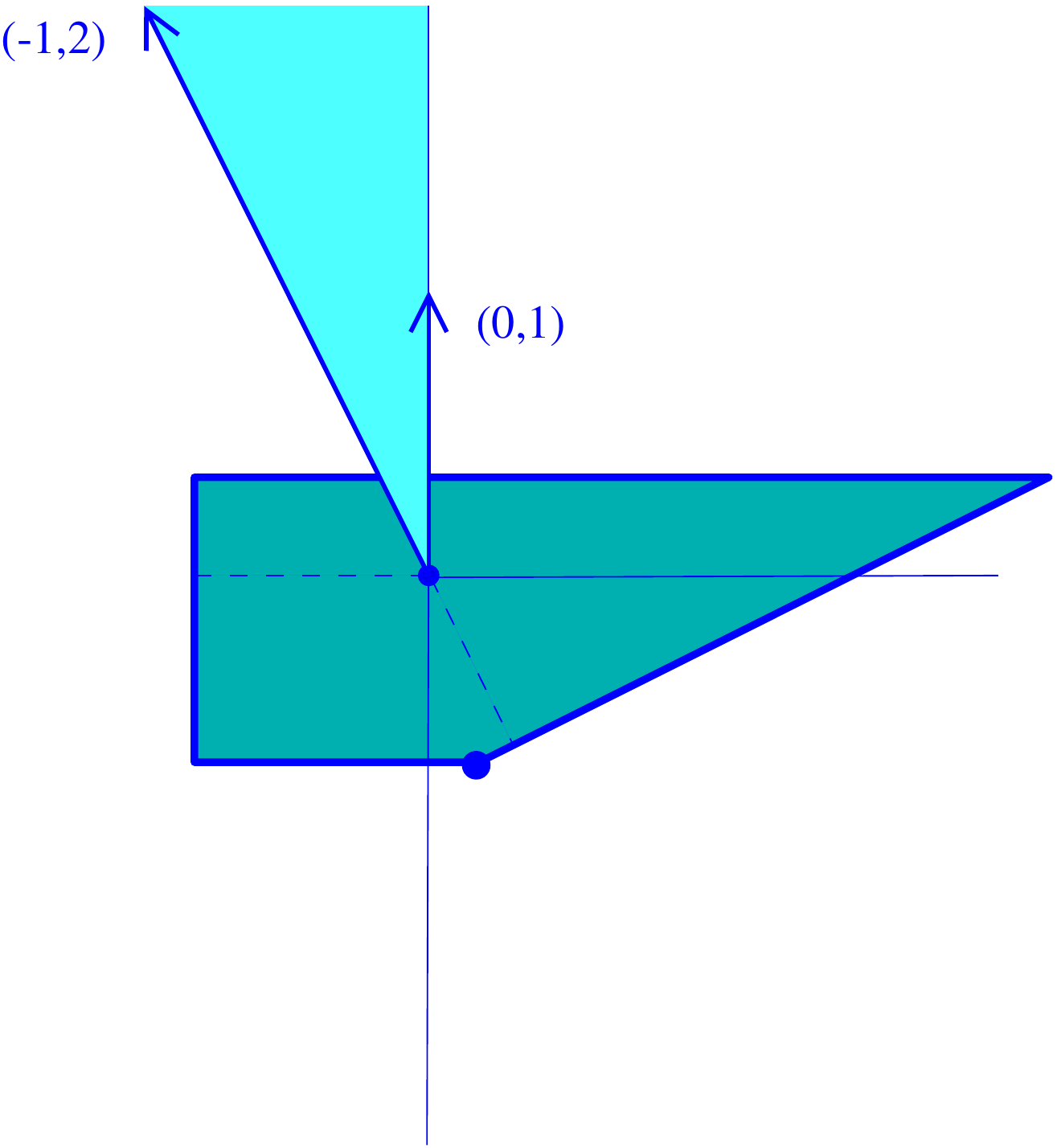}
\caption{A vertex of $P_2$ and its corresponding maximal cone.}
\label{smoothvertex}
\end{center}
\end{figure}

Hirzebruch surfaces were introduced by Hirzebruch in his thesis \cite{hirzebruch} and turn out in a number of different
contexts. They are complex algebraic surfaces, parametrized by positive integers. For each such integer $n$, the Hirzebruch surface
$\F_n$ is the projectivization of the bundle $\mathcal{O}\oplus\mathcal{O}(n)$ over $\C \P^1$. 
\begin{figure}[h]
\begin{center}
\includegraphics[width=5cm]{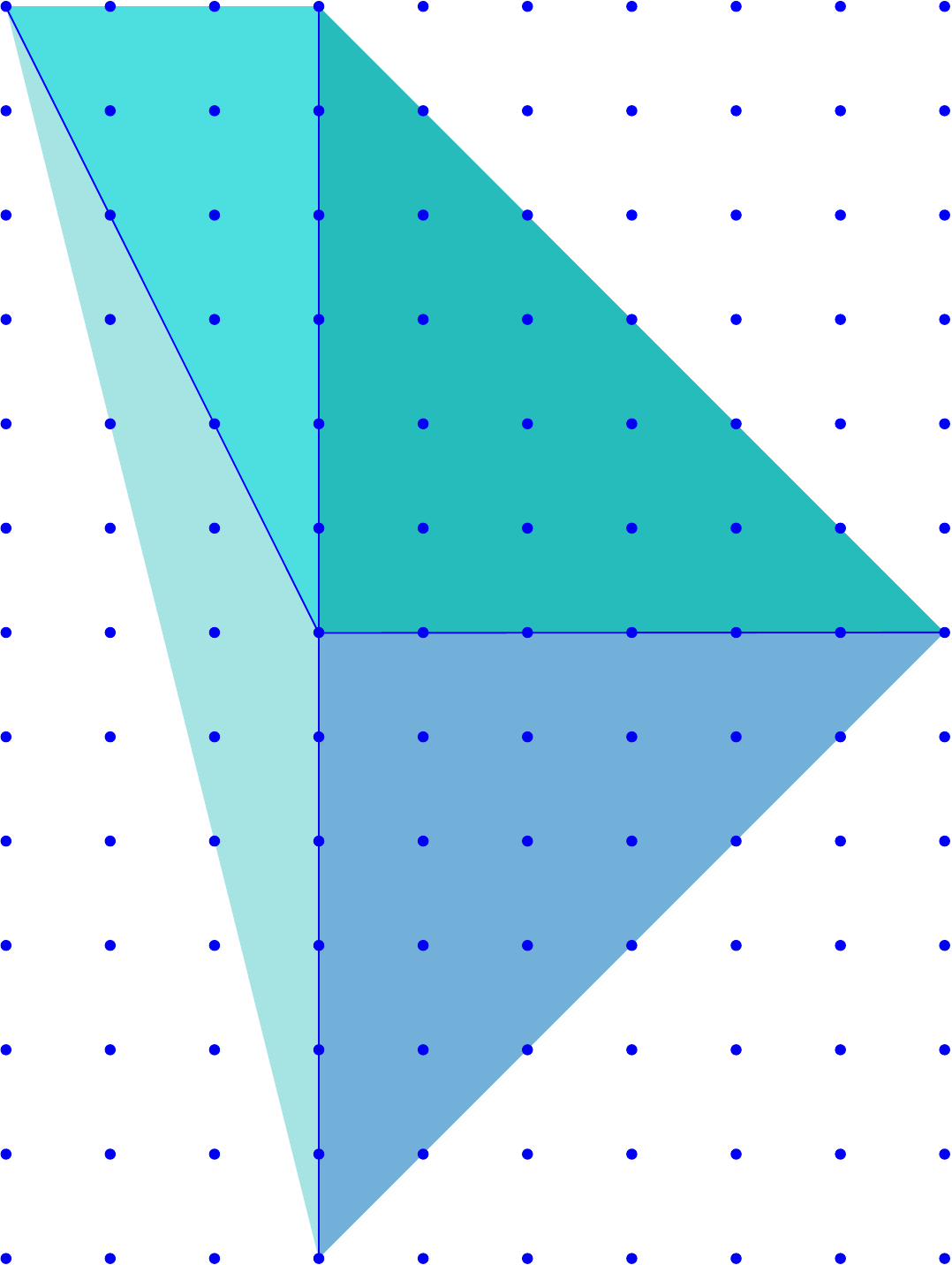}
\caption{The fan $\Delta_2$.}
\label{smoothfan}
\end{center}
\end{figure} 
Hirzebruch surfaces are toric manifolds. In fact, for each $n$, the surface $\F_n$ can be constructed from the fan 
$\Delta_n$ drawn in Figure~\ref{smoothfan} for $n=2$; the slanting ray passes through the point $(-1,n)$, 
while the other three are fixed (see, for example, \cite[Section~1.1]{fulton}).
The fan $\Delta_n$ is smooth. In fact,  $\Delta_n$ is rational, since each of its rays intersects the lattice $\Z^2$, and clearly simplicial; moreover, the primitive generators of the two rays of each of its maximal cones form a basis of $\Z^2$.
The fan $\Delta_n$ is also polytopal, i.e. there exists a convex polytope $P_n\subset(\R^2)^*$, 
given in this case by a right trapezoid, having $\Delta_n$ as normal fan.
Each ray of $\Delta_n$ is the inward--pointing normal ray of a facet of $P_n$ and the 
maximal cones  of $\Delta_n$ are in duality with the vertices of $P_n$ (see Figure~\ref{smoothpolytope}).
Since the normal fan to $P_n$ is smooth, the trapezoid $P_n$ is also smooth (see Figure~\ref{smoothvertex}).
We remark that there are infinitely many such trapezoids $P_n$.
The choice of a particular one defines a K\"ahler structure on $\F_n$. One way to see this is to recall that, from each $P_n$,  
one obtains  a symplectic toric manifold via the Delzant construction \cite{delzant}; this is equivariantly diffeomorphic 
to the toric manifold $\F_n$, and its symplectic structure is compatible with the complex one \cite[Chapter VI]{audin}. 
From now on, we will make a choice and call $P_n$ the right trapezoid of vertices $(0,0)$, $(0,1)$, $(n+1,1)$, and $(1,0)$. 

Classically, the above constructions make sense only for positive integers $n$. Is there a way to extend the Hirzebruch family, 
allowing $n$ to be any positive real number $a$?
 
\begin{figure}[h]
\begin{center}
\includegraphics[width=5cm]{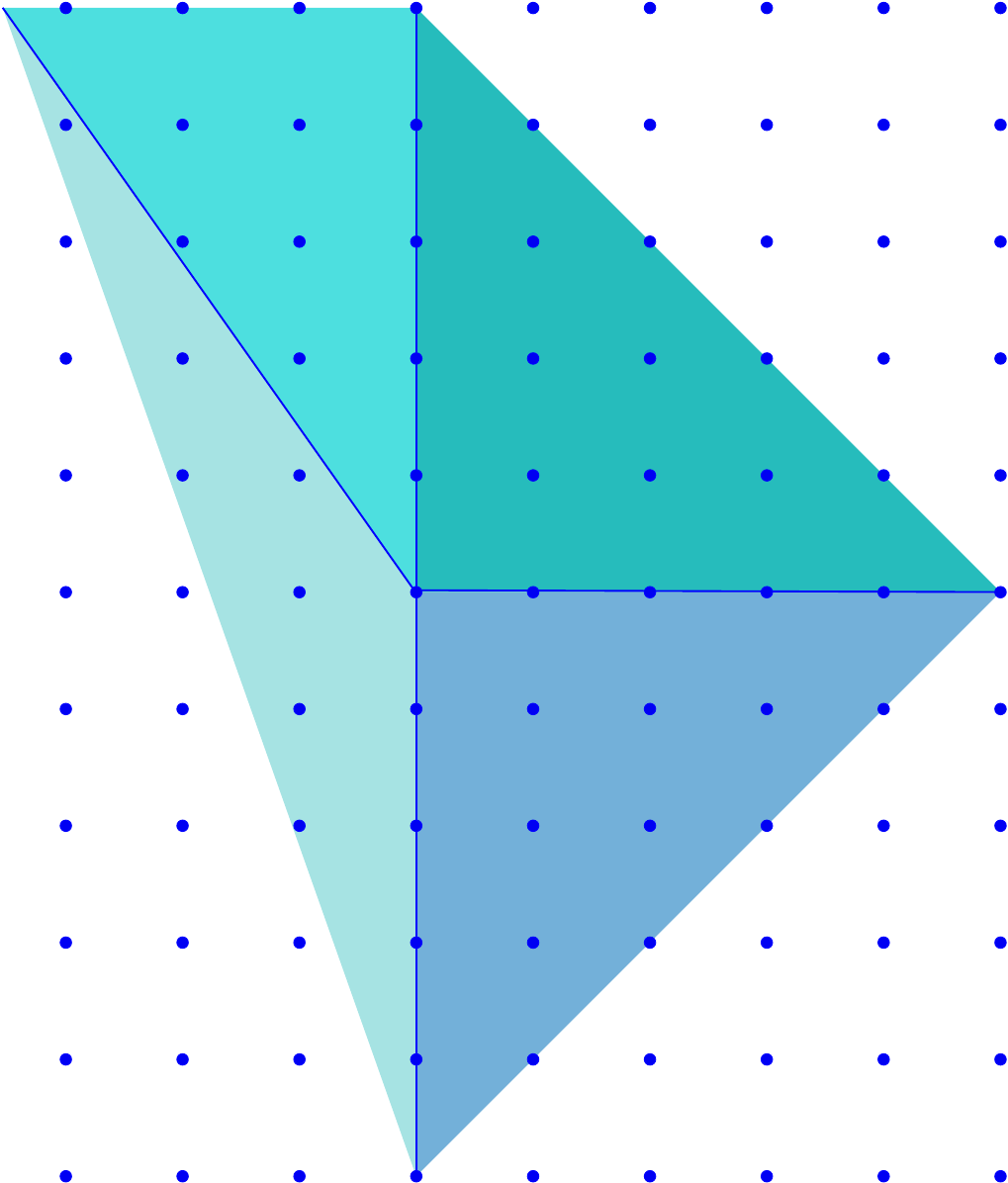}
\caption{The fan $\Delta_{\sqrt{2}}$.}
\label{quasihirzefan}
\end{center}
\end{figure}

From the convex--geometric viewpoint, for any positive real number $a$, we still have a simplicial polytopal fan $\Delta_a$, 
together with the corresponding right trapezoid $P_a$. However, for nonrational values of $a$, $\Delta_a$, and thus $P_a$, are not
rational in $\Z^2$ (see Figure~\ref{quasihirzefan}).
According to \cite{p,cx}, we may consider $\Delta_a$ and $P_a$ in a nonrational setting and construct, 
for each positive real number $a$, a  K\"ahler toric quasifold $\F_a$ (see Section~\ref{preliminaries}). 
When $a$ is rational, we obtain a K\"ahler orbifold; when $a$ is an integer $n$, we recover the Hirzebruch surface $\F_n$.
In this article, we  present two alternative constructions of $\F_a$. All three points of view rely on the idea, introduced in \cite{p}, 
of replacing the lattice $\Z^2$ by the quasilattice $$Q_a=\text{span}_\Z\{(1,0),(0,1),(0,-1),(-1,a)\}=\Z\times (\Z+ a\Z).$$

The first construction, introduced in \cite{bz1}, provides a smooth model for each $\F_a$. More precisely, for each $a$, 
we construct a compact, complex manifold $N_a$, endowed with a holomorphic, transversely K\"ahler foliation $\calF_a$, 
such that the leaf space is exactly $\F_a$. The general construction builds on previous work by Meersseman--Verjovsky; in particular, the case $a$ positive integer is treated in \cite[Example~5.6]{MV}. The manifolds $N_a$
are in the so--called LVM family: a large class of compact, complex, non--K\"ahler manifolds, introduced by Lopez de Medrano, Verjovsky, and Meersseman in \cite{ldm,M}. As shown in \cite{LN,M}, they are endowed with a holomorphic foliation.
The foliated manifolds $(N_a,\calF_a)$ that we obtain are of complex dimension $3$, with one--dimensional foliation.
The topological type of the generic leaf is  
$$\calF_{[\vz]}\simeq\left\{\begin{array}{ll}(S^1)^2 \;&\text{if}\,a\in\Q\\
\;S^1\times\R\;&\text{if}\,a\notin\Q\end{array}\right.;$$
it varies from closed to nonclosed, depending on whether $a$ is rational or not. When $a$ is irrational,
the closure of the generic leaf is diffeomorphic to $(S^1)^3$. This model has a symplectic counterpart, 
where $N_a$ is seen as a presymplectic foliated manifold \cite{bz1,rome}.
Finally, we remark that, for each $a$, it is possible to construct infinitely many further pairs $(N,\calF)$, 
where the manifold $N$ can have any dimension greater or equal to $3$, and such that the leaf space is $\F_a$. For a study of these manifolds see \cite{bz2}.
 
\begin{figure}[h]
\begin{center}
\includegraphics[width=8cm]{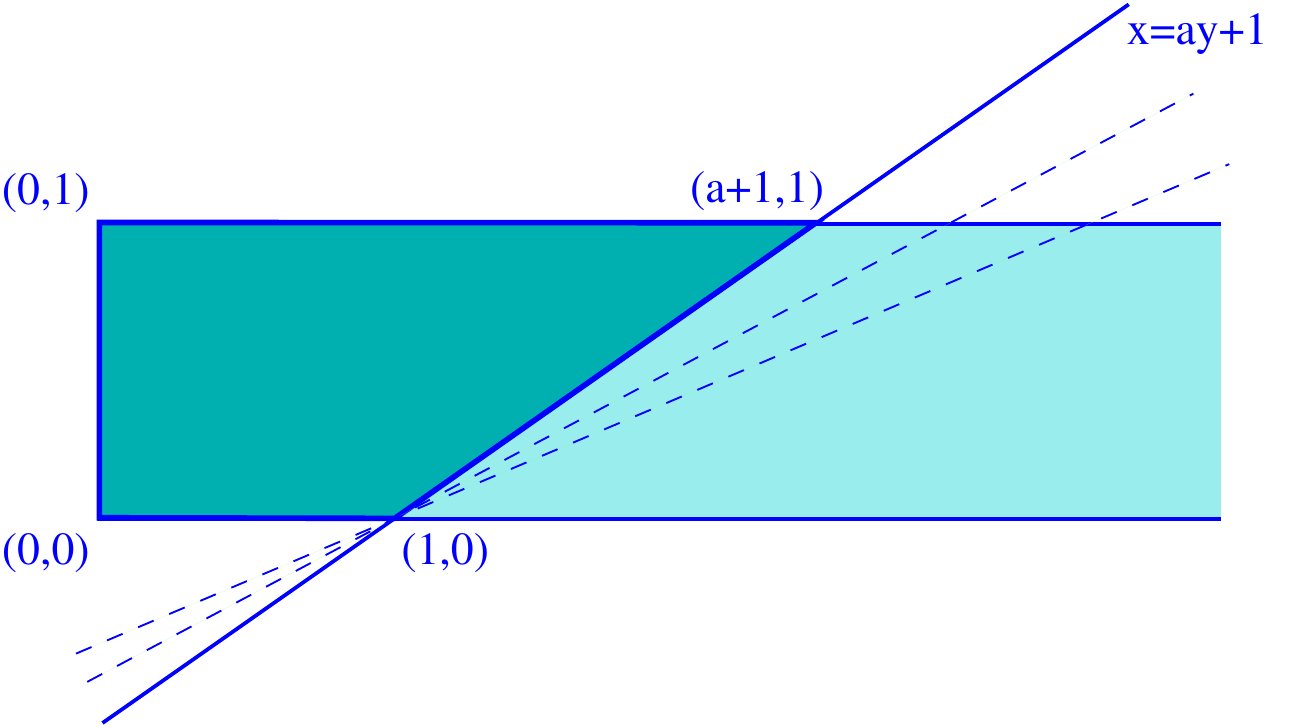}
\caption{Cutting $\C\times S^2$ in arbitrary directions.}
\label{hirzecut}
\end{center}
\end{figure}

The second construction consists in realizing the quasifold $\F_a$ as a symplectic cut of the manifold $\C\times S^2$. 
The idea is that the trapezoid $P_a$ can be obtained by cutting the strip
$[0,+\infty)\times [0,1]$ with the line $x=ay+1$ (see Figure~\ref{hirzecut}), where the strip indeed corresponds
to the symplectic toric manifold $\C\times S^2$.
However, for nonrational values of $a$, this line is irrationally sloped and the classical Lerman cut \cite{lerman} 
cannot be applied. We rely on a generalization of this procedure that allows cutting in an arbitrary direction \cite{cut}. 
One of the relevant features of this point of view is that it allows to express the symplectic quasifold 
$\F_a$ as a disjoint union of a compact two--dimensional quasifold and of an open dense subset. The first is given by $S^2/\Gamma_a$, where the countable group $\Gamma_a=Q_a/\Z^2$ acts on $S^2$ by rotations, 
around the $z$--axis, of integer multiples of $2\pi a$; the second is an open subset of $\C\times S^2$, 
again modulo the above action of $\Gamma_a$ on $S^2$.

As we have seen, the topological structure of $\F_a$ varies depending on the rationality of $a$. 
Notice that, within the family $\F_a$, we can pass from one Hirzebruch surface to the other;
however, this is not done via a deformation, in accordance with the fact that Hirzebruch surfaces $\F_n$ have two distinct diffeomorphism types, 
depending on the parity of $n$ \cite{hirzebruch}.

\section{Preliminaries}\label{preliminaries}
We recall the basic facts on toric geometry for nonrational convex polytopes, following \cite{p,cx}, and we apply them to the 
generalized Hirzebruch setting. Section~\ref{whatis} is an outline of some relevant notions in convex geometry. 
In Section~\ref{toricquasifolds}, we recall what {\em toric quasifolds} are and we describe the generalized Hirzebruch quasifold $\F_a$.

\subsection{Convex polytopes, fans, and quasilattices}\label{whatis}
A convex polytope $P\subset(\R^k)^*$ determines a fan $\Delta\subset\R^k$, known as the {\em normal fan} to $P$. 
There is an inclusion--reversing bijection between cones in $\Delta$ and faces of $P$.
For example, each of the {\em rays}, namely the one--dimensional cones, of $\Delta$ is orthogonal to a facet of $P$ and points towards its interior. 
The maximal cones of $\Delta$, on the other hand, correspond to the vertices of $P$.
For more on convex polytopes and their normal fans, we refer the reader to \cite{ziegler}. We just recall a few more relevant facts.
The convex polytope $P$ is said to be {\em simple} if its normal fan is {\em simplicial}, 
namely if each cone of $\Delta$ is a cone over a simplex. 
Therefore, a $k$--dimensional convex polytope $P$ is simple if, and only if, each of its vertices is the intersection of exactly $k$ facets. 
The convex polytope $P$ is said to be {\em rational} if its normal fan is, namely if there is a lattice $L\subset\R^k$ 
that has nonempty intersection with each of the rays of $\Delta$. 
Finally, the convex polytope $P$ is said to be {\em smooth} if its normal fan is, namely if 
$\Delta$ is rational, simplicial and, for each of the maximal cones of $\Delta$,
the primitive generators of its rays form a basis of the lattice. 
The toric variety corresponding to a smooth fan is also smooth.
The lattice and the set of primitive (or minimal) generators of the rays of any rational fan
are key ingredients in the construction of the corresponding toric variety (see, for example, \cite[Chapter 6]{audin}) or
\cite{cox}).
What happens for nonrational fans? 
The idea, introduced by the second author in \cite{p}, is to replace the lattice with a {\em quasilattice}
$Q\subset\R^k$, namely the $\Z$--span of a set of $\R$--spanning vectors of $\R^k$. 
One way of constructing a quasilattice in this setting is to choose a set of generators of the rays of the fan and to let 
$Q$ be equal to their $\Z$--span. There is of course a lot of freedom in the choice of the generators in general, but in some cases 
one can be guided by the underlying geometric setup. This is exactly what happens in the Hirzebruch setting.
The set of primitive generators for the rays of the Hirzebruch fan $\Delta_n$ is given by $\{(1,0),(0,1),(0,-1),(-1,n)\}$.
Therefore, if we allow $n$ to be any positive real number $a$, the natural choice for a set of generators for the rays of the fan $\Delta_a$
is $\{(1,0),(0,1),(0,-1),(-1,a)\}$. This yields the quasilattice $$Q_a=\text{span}_\Z\{(1,0),(0,1),(-1,a)\}=\Z\times (\Z+ a\Z)\supseteq \Z^2.$$ Notice that,
for $a$ rational, $Q_a$ is a lattice; for $a$ natural this lattice equals $\Z^2$.
\subsection{Toric quasifolds}\label{toricquasifolds}
Toric quasifolds are the natural generalization of toric manifolds in the nonrational setting. Quasifolds 
were introduced by the second author in \cite{p}. They are singular spaces that generalize manifolds and orbifolds. 
They are locally modeled by quotients of manifolds modulo the smooth action of countable groups and they are typically not Hausdorff. 
Notable examples of quasifolds are the so--called {\em quasitori} $D^k=\R^k/Q$, $Q$ being a quasilattice in $\R^k$, 
and its complexification $D^k_{\C}=\C^k/Q$.

We now recall the extensions to the nonrational setting of the Delzant construction \cite{delzant} and of its complex counterpart \cite{audin}.
Let $P\subset(\R^k)^*$ be a simple, dimension $k$ convex polytope, let $\{X_1,\ldots,X_d\}$ be a set of generators of the rays of its 
normal fan and let $Q$ be a quasilattice that contains them. According to \cite[Theorem 3.3]{p}, for each triple  $(P,\{X_1,\ldots,X_d\},Q)$ 
one can construct a symplectic, dimension $2k$ quasifold $M_P$ that is endowed with an effective Hamiltonian action of 
the quasitorus $D^k$; this action is Hamiltonian and the image of the corresponding moment mapping is the polytope $P$. 
On the other hand, by \cite[Theorem 2.2]{cx}, for the same triple one can also construct a complex, dimension $k$ quasifold $X_P$ 
that is endowed with a holomorphic action of the  complex quasitorus $D^k_{\C}$; this action has a dense open orbit. 
Finally, by \cite[Theorem 3.2]{cx}, $M_P$ and $X_P$ are equivariantly diffeomorphic and the induced symplectic form on $X_P$ is K\"ahler. 
The spaces $M_P$ and $X_P$ are respectively called the {\em symplectic} and {\em complex toric quasifold} corresponding 
to the triple $(P,\{X_1,\ldots,X_d\},Q)$.

\begin{figure}[h]
\begin{center}
\includegraphics[width=7cm]{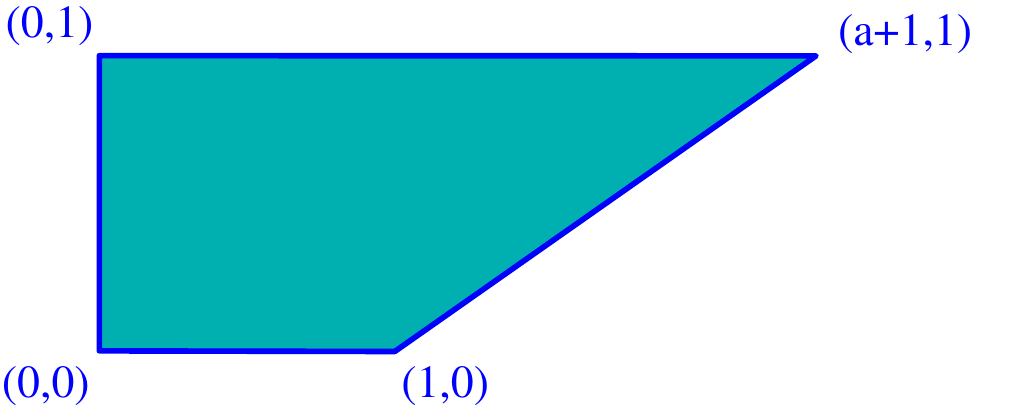}
\caption{The trapezoid $P_a$.}
\label{quasihirze}
\end{center}
\end{figure}

Let us now consider, for any positive real number $a$, the generalized Hirzebruch trapezoid $P_a$; its vertices are given by 
$(0,0)$, $(0,1)$, $(a+1,1)$, and $(1,0)$ (see Figure~\ref{quasihirze}). 
Take the countable group $$\Gamma_a=Q_a/\Z^2\simeq (\Z+a\Z)/\Z.$$
One can verify that, for the triple 
$$\Big(P_a,\{(1,0),(0,1),(0,-1),(-1,a)\}, Q_a\Big),$$
the results above yield the symplectic toric quasifold
$$\frac{\{\, \vz\in\C^4\;|\; |z_1|^2+a|z_3|^2+|z_4|^2=1+a, |z_2|^2+|z_3|^2=1\,\}}{\{\,(e^{2\pi i r}, e^{2\pi i s}, e^{2\pi i (s+ar)}, e^{2\pi i r})\in \R^4/\Z^4\;|\;r,s\in \R\,\}},$$
acted on by the quasitorus
$D^2_a=\R^2/Q_a\simeq S^1\times (S^1/\Gamma_a)$, and the complex toric quasifold
$$\frac{\{\vz\;|\;z_3z_4\neq0\}\cup\{\vz\;|\;z_1z_3\neq0\}\cup
\{\vz\;|\;z_1z_2\neq0\}\cup\{\vz\;|\;z_2z_4\neq0\}}{\{\,(e^{2\pi i u}, e^{2\pi i v}, e^{2\pi i (v+au)}, e^{2\pi i u})\in \C^4/\Z^4\;|\;u,v\in \C\,\}},$$
acted on by the complex quasitorus $(D^2_a)_{\C}=\C^2/Q_a\simeq \C^*\times (\C^*/\Gamma_a)$.
We denote $\F_a$ the resulting K\"ahler toric quasifold. If $a$ is a rational number that is not an integer, we can write
$a=\frac{p}{q}$, with $p,q$ coprime positive integers and $q>1$. In this case,
$\F_a$ is a topologically smooth K\"ahler toric orbifold with two disjoint singular divisors of order $q$ given by $z_2=0$ and $z_3=0$, which correspond
to the bases of the trapezoid $P_a$.
Notice finally that, for $a$ equal to a positive integer $n$, we obtain the standard Hirzebruch surface $\F_n$.
\section{Foliations modeling $\F_a$} 
In this section, we apply a construction developed in \cite{bz1} that allows to model complex and symplectic 
quasifolds by complex and presymplectic, foliated, smooth manifolds, respectively. This viewpoint builds 
on the articles \cite{p} and \cite{MV}. In Section~\ref{vectorconfig}, convex data are interpreted in the context
of vector configurations. In Section~\ref{complex}, we review the complex construction focussing on the Hirzebuch
family and finally, in Section~\ref{symplectic}, we illustrate the symplectic side of the picture, for which we refer also to \cite{rome}.

\subsection{Vector and point configurations}\label{vectorconfig}
Consider a complete simplicial fan $\Delta\subset\R^k$; complete means that the union of its cones is $\R^k$.
Let $\{X_1,\ldots,X_d\}$ be a choice of generators of the rays of $\Delta$ and let $Q$ be a quasilattice containing them.
Each triple $(\Delta,\{X_1,\ldots,X_d\},Q)$ can be encoded in a well--studied convex object, a  
{\em triangulated vector configuration} \cite[Section~2.1]{bz1}. In the Hirzebruch case, the triple $(\Delta_n,\{(1,0),(0,1),(0,-1),(-1,n)\}, Q_n)$ can be encoded in $(V'_n,\bT)$, where
$V'_n=((1,0),(0,1),(0,-1),(-1,n))$
and
$\bT=\{\{1,2\},\{2,4\},\{3,4\},\{1,3\},\{1\},\{2\},\{3\},\{4\},\varnothing\}.$
The set of vectors $V'_n$ is a {\em vector configuration:} a finite, ordered list of vectors, allowing repetitions.
The vector configuration $V'_n$ contains the following information: a set of ray generators, and therefore the rays themselves, 
and the lattice $\text{span}_{\Z}\{V'_n\}=\Z^2$. The set $\bT$ is a {\em triangulation} of $V'_n$. It is a collection of subsets of $\{1,2,3,4\}$, with suitable properties. In this case, $\bT$ carries the combinatorial information of $\Delta_n$: for example, $\{1,2\}$ corresponds to the maximal cone in $\Delta_n$
generated by the first and second vectors of $V_n'$. By convention, the cone corresponding to $\varnothing$ is $\{0\}$.

In general, for our construction, we will need a vector configuration that is {\em balanced}, meaning that the 
sum of its vectors is $0$, and {\em odd}, meaning that $\text{card}(V'_n)-\text{dim}(\text{span}_{\R}(V'_n))$ is odd.
We can ensure that these two assumptions are verified by adding to our configuration, if necessary, some extra vectors,
known as {\em ghost vectors}. For the Hirzebruch fan we can take 
$$V_n=\Big((1,0),(0,1),(0,-1),(-1,n),(0,-n)\Big).$$
The fifth vector $(0,-n)$ is a ghost vector since is not indicized by any subset of $\bT$. 
Now $\text{card}(V_n)-\text{dim}(\text{span}_{\R}(V_n))=5-2=3=2m+1$, with $m=1$. 
Remark that the configuration $V_n$ is not uniquely determined. For example, we could append to $V_n$ any even number of ghost vectors in $\Z^2$. This operation is sometimes necessary in order to ensure that the $\Z$--span of the vector configuration is the given lattice or quasilattice.

\begin{figure}[h!]
\begin{center}
\includegraphics[width=7cm]{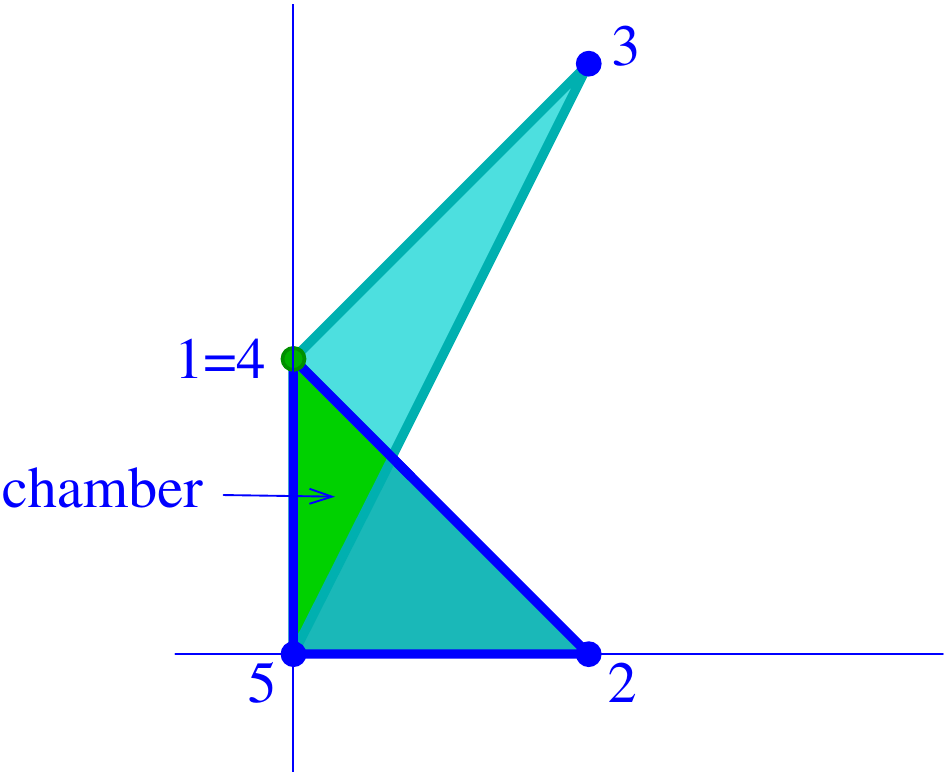}
\caption{Point configuration and chamber.}
\label{figurechamber}
\end{center}
\end{figure}

We consider, as initial convex datum, $(V_n,\bT)$, with $n$ positive integer.  Its natural generalization to positive real numbers $a$ is 
the triangulated vector configuration 
$(V_a,\bT),$ with $V_a=((1,0),(0,1),(0,-1),(-1,a),(0,-a)).$ 
Now apply {\em Gale duality}: choose a basis of the relations $\text{Rel}(V_a)\subset\R^5$ of $V_a$ and write it as rows in a matrix:
$$\left(\begin{array}{ccccc}
1&1&1&1&1\\
0&1&1&0&0\\
1&0&a&1&0
\end{array}\right).
$$
Ignoring the first row, we interpret each column as the real and imaginary parts of a complex number. 
This yields the {\em configuration of points} (that is a finite ordered list)
$\Lambda=(i,1,1+ia,i,0)$
in affine space $\C^m=\C$, as shown in Figure \ref{figurechamber}.
The corresponding combinatorial datum is a {\em virtual chamber} $\bT^*$, namely a collection of subsets of 
$\{1,2,3,4,5\}$ satisfying certain properties. In our case, we obtain:
$$\bT^*=\{\{3,4,5\},\{1,3,5\},\{1,2,5\},\{2,4,5\}\}.$$
In fact, the intersection of the four triangles determined by $\bT^*$ is nonempty. This happens if, and only if, the fan is polytopal \cite{S}.

As we have seen, there are infinitely many triangulated vector configurations encoding a given triple. Moreover,
given $(V,\bT^*)$, its Gale dual $\Lambda$ is not unique, rather it is determined up to a real affine automorphism of 
$\C$ (see \cite[Section~2.2.2]{bz1} and also \cite[Section~ 1.2~]{rome} for details).
For an exhaustive treatment of notions like vector and point configurations, Gale duality, triangulations, and chambers, 
we refer the reader to \cite{DL-R-S}. 
\subsection{The complex foliated manifolds $(N_a,\calF_a)$}\label{complex}
We consider the LVM datum $(\Lambda,\bT^*)$.
The chamber $\bT^*$ determines the open subset $U(\bT^*)$ of $\C \P^4$ given by the projectivization of the open subset of $\C^5$ given by:
$$\{\vz\;|\;z_3z_4z_5\neq0\}\cup\{\vz\;|\;z_1z_3z_5\neq0\}\cup
\{\vz\;|\;z_1z_2z_5\neq0\}\cup\{\vz\;|\;z_2z_4z_5\neq0\}.
$$
The configuration $\Lambda$ determines the following $\C_{\Lambda}$--action on $U(\bT)$:
$$\begin{array}{ccc}
\C_{\Lambda}\times U(\bT^*)&\longrightarrow&U(\bT^*)\\
(t,[ z_1:z_2:z_3:z_4:z_5 ])&\longmapsto&[ e^{-2\pi t}z_1:e^{2\pi i t}z_2:e^{2\pi i (1+ia)t}z_3:e^{-2\pi t}z_4:z_5 ].
\end{array}
$$
The quotient $N_a=U(\bT^*)/\C_{\Lambda}$ is a compact, complex manifold  \cite[Section~2.2.4]{bz1}. 
Remark that the above procedure applies to nonpolytopal fans as well \cite{bz1}, yielding a generalized LVM manifold \cite{bosio}.
Following \cite{U}, consider now the conjugate point configuration, $\overline{\Lambda}$, of $\Lambda$.
The action of $\C_{\overline{\Lambda}}$ on $U(\bT^*)$ commutes with that of $\C_{\Lambda}$; therefore, there is an induced
action on $N_a$: 
$$\begin{array}{ccc}
\C_{\overline{\Lambda}}\times N_a&\longrightarrow& N_a\\
(t,[ z_1:z_2:z_3:z_4:z_5 ])&\longmapsto&[ e^{2\pi t}z_1:e^{2\pi i t}z_2:e^{2\pi i (1-ia)t}z_3:e^{2\pi t}z_4:z_5 ].
\end{array}
$$
Its orbits give rise to a holomorphic foliation $\calF_a$ in $N_a$ \cite[Section~2.3]{bz1}. The set of generic points in 
$N_a$ is the dense, open orbit of the induced action of $(\C^*)^4$ on $N_a$. Moreover, $\C_{\overline{\Lambda}}$ 
acts on $N_a$ as a subgroup of $(S^1)^4\subset(\C^*)^4$. This implies that the foliation $\calF_a$ is Riemannian \cite[Section~2.3.2]{bz1}.

The topological type of the generic leaf and of its closure depends on the measure of the rationality of $V$ \cite[Section~2.1.1]{bz1}. 
At a generic point $[\vz]\in N_a$, we obtain:
$$\calF_{[\vz]}\simeq\left\{\begin{array}{ll}(S^1)^2 \;&\text{if}\,a\in\Q\\
\;S^1\times\R\;&\text{if}\,a\notin\Q\end{array}\right.$$
$$\overline{\calF}_{[\vz]}\simeq\left\{\begin{array}{ll}(S^1)^2\;&\text{if}\,a\in\Q\\
(S^1)^3\;&\text{if}\,a\notin\Q.\end{array}\right.$$
This shows very clearly how differently leaves behave when $a$ passes from rational to irrational numbers. 
This, of course, reflects on the topology of the leaf space $\F_a$. In our case,  the holomorphic projection from $(N_a,\calF)$ 
to the complex leaf space $\F_a$ has a very simple expression:
$$\begin{array}{ccc}N_a&\longrightarrow& \F_a\\
\left[ z_1:z_2:z_3:z_4:z_5 \right]&\longmapsto&\left[ z_1z^{-1}_5:z_2z^{-1}_5:z_3z^{-1}_5:z_4z^{-1}_5 \right].
\end{array}
$$
We recall that the complex structure of the leaf space does not depend on the choice of Gale dual. More generally, we have seen that 
there are infinitely many triangulated vector configurations $(V,\bT)$ encoding the triple 
$$\Big(\Delta_a,\{(1,0),(0,1),(0,-1),(-1,a)\},Q_a\Big).$$  For each of them, there are infinitely many choices 
of Gale dual point configurations $(\Lambda,\bT^*)$. Thus, we obtain infinitely many foliated manifolds, but only one complex leaf space \cite[Theorem~2.1]{rome}:
$$\xymatrix{
(V,\bT)\ar@<1ex>[d]\ar@{-->}[r]&(\Lambda,\bT^*)\ar@<1ex>[l]\ar[r]&N/\calF\ar[d]\\
(\Delta_a,\{(1,0),(0,1),(0,-1),(-1,a)\},Q_a)\ar@{-->}[u]\ar[rr]&&\F_a
}
$$
Dashed arrows here represent directions in which we make choices. On the other hand, the complex structure of the leaves 
does vary upon the choice of Gale dual. If $a$ is irrational, the leaf $\calF_{[\vz]}$ is $\C^*$ for $z_2z_3\neq0$; it is a compact complex torus otherwise. If $a$ equals $\frac{p}{q}$, with $p,q$ coprime positive integers, 
$$\calF_{[\vz]}=\left\{\begin{array}{lcl}
\C/(\Z\oplus iq\Z) & \text{for} & z_2z_3\neq0\\
\C/\Z\oplus i\Z & \text{for} & z_2=0\;\text{or}\; z_3=0.
\end{array}
\right.
$$
By varying our choice of Gale dual, we obtain all possible two--dimensional compact complex tori, in accordance with \cite[Theorem~G]{MV}. 
Finally, the generic leaf is a $q$--sheeted cover of the leaf $\calF_{[\vz]}$ through 
$[\vz]=[z_1:0:z_2:z_3:z_4:z_5]$: when approaching the point $[\vz]$, the generic leaf winds $q$ times around the leaf 
$\calF_{[\vz]}$, compatibly with \cite[Corollary~B]{MV}. According to \cite{klmv}, the leaf space $\F_a$, obtained from a nonrational LVM datum $(\Lambda,\bT^*)$, can be seen from the viewpoint of noncommutative geometry.
\subsection{The presymplectic foliated manifolds}\label{symplectic}
Let $(\Delta,\{X_1,\ldots,X_d\},Q)$ be a triple such that $\Delta$ is the normal fan
of a simple, convex polytope $P$. Let $(V,\bT)$ be a triangulated vector configuration encoding the given triple.
The resulting foliated manifold $(N,\calF)$ can be viewed in a symplectic setting, by applying a simple variant of the generalized Delzant construction \cite[Theorem~3.3]{p} 
introduced in \cite[Proposition~4.3]{rome}.
Focussing on our Hirzebruch family, consider the triple $(P_a,\{(1,0),(0,1),(0,-1),(-1,a)\}, Q_a)$, 
with the additional datum of the half--plane $-ay\geq -2a$.
We then obtain the connected subgroup of $(S^1)^5$ given by $\exp(\text{Rel}(V_a))$; its induced action on $\C^5$ is Hamiltonian, with moment mapping 
$$\Psi_a(\vz)=\Big(|\vz|^2-2(a+1), |z_1|^2+a|z_3|^2+|z_4|^2-1+a, |z_2|^2+|z_3|^2-1 \Big).$$
We find the following diagram:
$$\xymatrix{
\Psi_a^{-1}(0)\ar[r]_{S^1}&\Psi_a^{-1}(0)/S^1\ar[d]_{\R^2}&\simeq&N_a\ar[d]^{\C_{\overline{\Lambda}_a}}\\
&(\F_a,\omega)&\simeq&\F_a}
$$
where $S^1$ acts diagonally on $\C^5$ and 
the $\R^2$--action on $\Psi_a^{-1}(0)/S^1$ is given by: $$(r,s)\cdot[\vz]=[e^{2\pi i r}z_1:e^{2\pi i s}z_2:
e^{2\pi i (s+ar)}z_3:e^{2\pi i r}z_4:z_5].$$ The manifolds $\Psi_a^{-1}(0)$ and $\Psi_a^{-1}(0)/S^1$ are both presymplectic.
The action of $\R^2$ on $\Psi^{-1}(0)/S^1$ is Hamiltonian and induces a foliation which is sent diffeomorphically
to $\calF_a$ \cite[Section~4]{rome}. The presymplectic structure of $\Psi_a^{-1}(0)/S^1$ defines a transversely
K\"ahler structure on $(N_a,\calF)$, a well known result \cite[Theorem~7]{M}. This presymplectic viewpoint 
was already investigated, and key, in \cite{MV}, in the rational setting. 
It is also related to recent articles by Lin--Sjamaar \cite{sjamaar} and Nguyen--Ratiu \cite{ratiu}. 
In particular, in the former, symplectic quasifolds are also viewed as leaf spaces of presymplectic manifolds.
\section{Arbitrary toric cuts of $\C\times S^2$}
It is well known that the Hirzebruch surface $\F_n$ can be obtained from weighted projective space $\C\P^2_{(1,n,1)}$ by 
blowing up its only singular point. In the symplectic category, this can be done by means of the symplectic blowing--up 
procedure. From the viewpoint of symplectic toric geometry, $\C\P^2_{(1,n,1)}$ corresponds to the triangle 
$T_n$ of vertices $(0,-1/n)$, $(0,1)$, and $(n+1,1)$, with the singular point mapping to the first vertex.
\begin{figure}[h]
\begin{center}
\includegraphics[width=7cm]{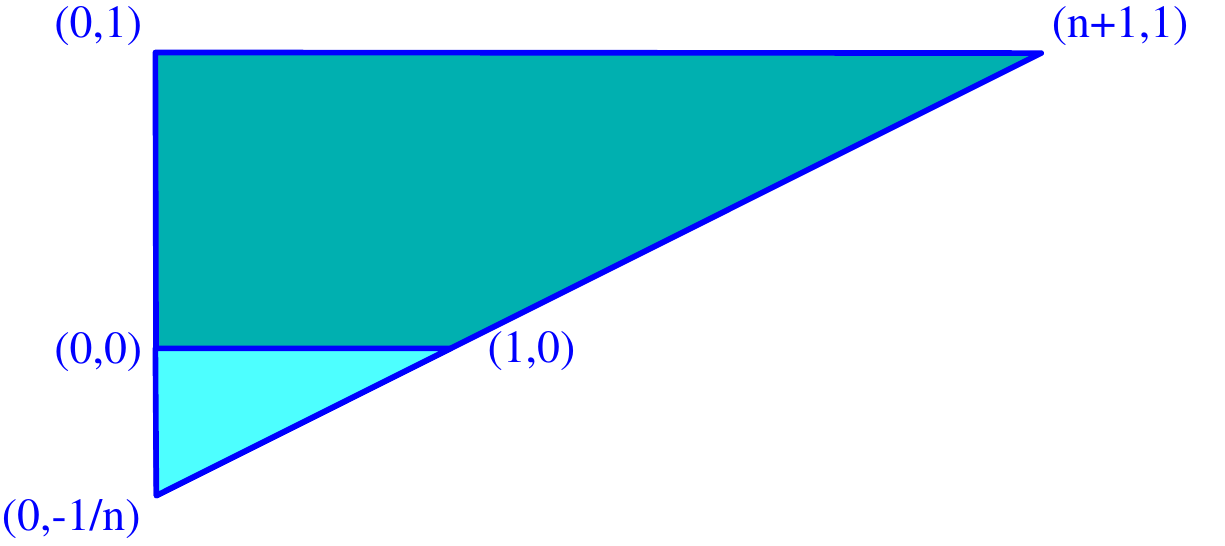}
\caption{Blowing--up $\C\P^2_{(1,n,1)}$.}
\label{smoothhirzetriangle}
\end{center}
\end{figure} 
 As it turns out, blowing--up this point of the amount $\frac{1}{n}$ corresponds to cutting off a corner of the triangle around the vertex, 
as in Figure~\ref{smoothhirzetriangle} (see, for example, \cite[Theorem 1.12, Example 1.22]{guilleminlibro}). 
The resulting convex polytope is the trapezoid $P_n$, which, as we know, corresponds to  $\F_n$. 
For an interesting application of this viewpoint, see \cite{gauduchon}.

\begin{figure}[h]
\begin{center}
\includegraphics[width=8cm]{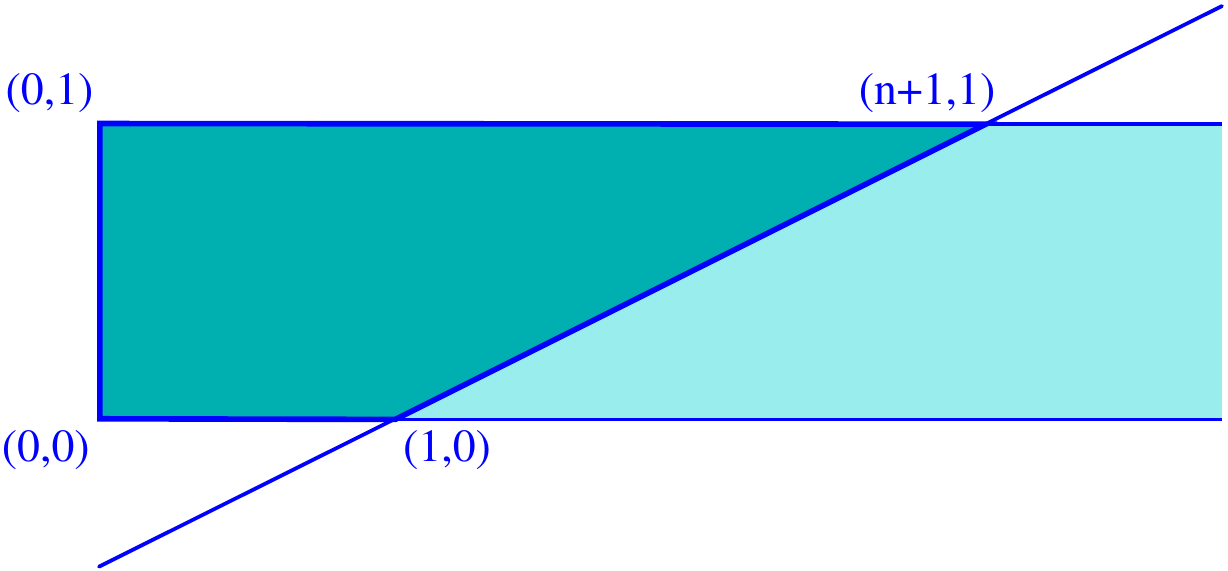}
\caption{Cutting $\C\times S^2$.}
\label{smoothhirzecut}
\end{center}
\end{figure}

Here we want to suggest an alternative approach to obtaining $\F_n$, using the symplectic cutting procedure. This operation was introduced 
by Lerman \cite{lerman}, and is an important generalization of the symplectic blowing--up operation. In the case of symplectic toric manifolds, 
the cutting procedure amounts to cutting the moment polytope with a hyperplane, and to considering the toric manifolds corresponding to 
the two cut polytopes -- this will only work if the latter are smooth. The basic idea in the Hirzebruch setting is that one can obtain the 
trapezoid $P_n$ by cutting the strip $[0,+\infty)\times [0,1]$ with the line $x=ny+1$, as in Figure~\ref{smoothhirzecut}.  
The set $[0,+\infty)\times [0,1]$ is not a polytope, however it is a pointed polyhedron, and the Delzant construction can 
still be applied \cite[Theorem 1.1]{cut}. In fact, from the symplectic viewpoint, $[0,+\infty)\times [0,1]$ corresponds 
to the noncompact toric manifold $\C\times S^2$ with its standard symplectic structure and torus action. 

Now, what happens if we replace $n$ with any positive real number $a$? In order to obtain the trapezoid $P_a$, we would need to cut the
strip $[0,+\infty)\times [0,1]$ with the line $x=ay+1$ (see Figure~\ref{hirzecut}). 
However, for $a$ irrational, the standard cutting procedure does not allow this. We apply a generalization of this 
operation to symplectic toric quasifolds \cite{cut}, which allows, in particular, to also make sense of cutting any 
symplectic toric manifold in an arbitrary direction. The main reason why we cannot cut our strip with the irrationally 
sloped line $x=ay+1$ is that the line does not have a normal vector sitting in $\Z^2$. So the idea is to add this vector and to replace
the lattice $\Z^2$ with the quasilattice 
$$\text{span}_\Z\{\Z^2\cup\{(-1,a)\}\}.$$
Notice that this equals once more the quasilattice $Q_a=\Z\times (\Z+a\Z)$. Following \cite[Section~5]{cut},
in order to cut the manifold $\C\times S^2$,  we first need to consider the symplectic toric quasifold corresponding to the triple 
$$\Big([0,\infty)\times [0,1], \{(1,0),(0,1),(0,-1)\},Q_a\Big).$$ 
One can verify directly, or deduce from \cite[Proposition 1.2]{cut}, that this is given by the quasifold $\C\times (S^2/\Gamma_a)$, 
where $\Gamma_a=Q_a/\Z^2\simeq\{\, e^{2\pi i al}\,|\, l\in \Z\,\}$ 
is acting on $S^2=\{(v,z)\in\C\times\R\;|\;|v|^2+z^2=1\}$ by 
rotations around the $z$--axis. The quasitorus acting on $\C\times (S^2/\Gamma_a)$ is the same $D^2_a\simeq S^1\times (S^1/\Gamma_a)$ of Section~\ref{toricquasifolds}.
We then apply the generalized cutting procedure of \cite[Section~2]{cut}, which in this case yields the following. 
First, we consider the Hamiltonian action of $S^1$ on $\C\times (S^2/\Gamma_a)$
given by $e^{2\pi i t}\cdot(u,[v:z])=(e^{-2\pi i t}u,[e^{2\pi i at}v:z])$. One verifies that its moment map is 
$$\Phi_Y(u,[v:z])=-|u|^2+\frac{a(z+1)}{2}.$$
Then we let $S^1$ act on $\C\times (S^2/\Gamma_a)\times \C$ with weight $-1$ on the last component; the corresponding moment mapping is given by 
$$\nu_{-}(u,[v:z],w)=-|u|^2+\frac{a(z+1)}{2}-|w|^2.$$ The cut space corresponding to the trapezoid is then given by the quotient
$\nu_{-}^{-1}(-1)/S^1$. By \cite[Theorem 2.3]{cut}, it coincides with the symplectic toric quasifold $\F_a$.
Moreover, by \cite[Remark 2.4]{cut}, this cut can be written, as in the classical cutting procedure, as the disjoint union of the 
quotient $\Phi_Y^{-1}(-1)/S^1$ and of the open dense subset
$$
\begin{array}{lll}
&\{\,(u,[v:z])\in \C\times (S^2/\Gamma_a)\;|\; \Phi_Y(u,[v:z])>-1\,\}=\\
&\{\,(u,[v:z])\in \C\times (S^2/\Gamma_a)\;|\; |u|^2-\frac{a(z+1)}{2}<1\,\}.
\end{array}
$$
We remark that the former is given by $S^2/\Gamma_a$ (see \cite[Example 3.2]{reduction}, while
the latter equals the open subset of $\C\times S^2$ given by
$$\{\,(u,v,z)\in \C\times S^2\;|\; |u|^2-\frac{a(z+1)}{2}<1\,\},$$
modulo the action of $\Gamma_a$. It is indeed a general fact (see \cite[Proposition~2.5]{cut}) that a dense, open subset of the cut space
is symplectomorphic to an open subset of the initial manifold, modulo the action of the countable group given 
by the quotient of the quasilattice by the initial lattice.
If $a$ is again equal to $\frac{p}{q}$, with $p,q$ coprime positive integers and $q>1$, the group $\Gamma_a$ equals $\Z/q\Z$.
It is interesting to notice that the orbifold $\F_a$ inherits its order $q$ singular sets from the orbifold $\C\times (S^2/(\Z/q\Z))$.

\begin{figure}[h!]
\begin{center}
\includegraphics[width=7cm]{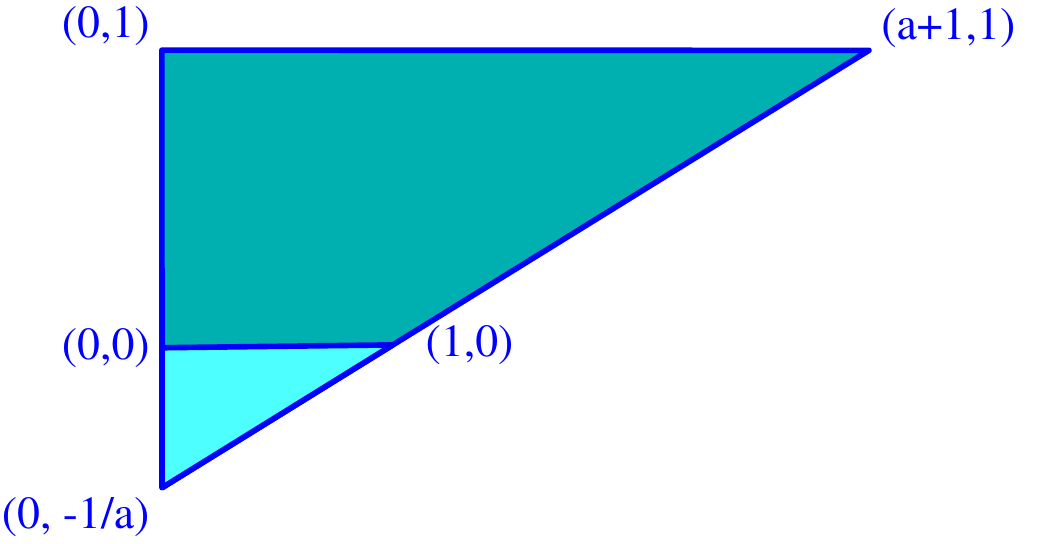}
\caption{Obtaining $\F_a$ via blowing--up.}
\label{hirzetriangle}
\end{center}
\end{figure}

We conclude by remarking that $\F_a$ can also be obtained as follows. 
First we consider the toric quasifold corresponding to the triple
$$\Big(T_a, \{(1,0),(0,-1),(-1,a)\},Q_a\Big),$$ where $T_a$ is the triangle of vertices $(0,-1/a)$, $(0,1)$ and $(a+1,1)$. 
Then, following \cite[Section~4]{cut}, we blow--up the point mapping to the first vertex of an amount $\frac{1}{a}$ (see Figure~\ref{hirzetriangle}). 
For $a$ equal to a positive integer $n$, this corresponds to obtaining the Hirzebruch surface by blowing--up weighted 
projective space, as explained at the beginning of this section.
\section*{Acknowledgements}
The research of the first two authors was partially supported by grant PRIN 2015A35N9B\_$\!$\_013 (MIUR, Italy) and by GNSAGA (INdAM, Italy).

\bigskip

{\small 

\noindent Fiammetta Battaglia, Elisa Prato\\
Dipartimento di Matematica e Informatica "U. Dini", Universit\`a di Firenze \\
Viale Morgagni 67/A, 50134 Firenze, ITALY

\noindent
fiammetta.battaglia@unifi.it, elisa.prato@unifi.it}

\medskip

{\small 

\noindent
Dan Zaffran\\ 
College of Marin\\ 
835 College Ave, Kentfield, CA 94904, USA

\noindent
dan.zaffran@gmail.com}

\end{document}